\def\BlackBox{\hspace*{\fill}\vrule width 4pt height 4pt depth 0pt}

\newcommand{\inprod}[2]{\left <#1,#2\right>}
\newcommand{\psd}[1]{\mathbb{S}^{#1}_+}

\newcommand{\zero}[0]{{\bf 0}}
\newcommand{\interior}[1]{\mbox{\bf int}\,#1}
\newcommand{\onesecond}[0]{\frac{1}{2}}
\newcommand{\tr}[1]{\mbox{\bf tr\,}#1}
\renewcommand{\det}[1]{\mbox{\bf det\,}#1}
\newcommand{\norm}[1]{\|#1\|}

\newcommand{\map}[0]{\longrightarrow}

\documentclass[runningheads]{llncs}
\usepackage[T1]{fontenc}
\usepackage{graphicx}
\usepackage[english]{babel}
\usepackage{epsf}
\usepackage{amssymb}
\usepackage{algorithmic}

\newtheorem{algorithm}{Algorithm}

\begin{document}
\title{A polynomial projective algorithm\\ for convex feasibility problems\\ with positive-definite constraints} 
\titlerunning{A projective algorithm for SDP}
%


\author{Sergei Chubanov\inst{1}\orcidID{0000-0003-3164-3415}}

\authorrunning{Sergei Chubanov}

\institute{Robert Bosch GmbH, Bosch Center for Artificial Intelligence, and University of Siegen, Germany\\ \email{sergei.chubanov@uni-siegen.de}}

%
\maketitle              
\begin{abstract}

We study a class of projective transformations of spectraplexes associated with self-dual cones and, on this basis, propose a polynomial-time algorithm for convex feasibility problems with positive definite constraints. At each iteration of the algorithm, either a feasible solution is found or a suitable valid inequality inducing a projective transformation allowing to bring the solution set closer to the center of an associated spectraplex. The closeness to the center is measured in terms of a potential function. The running time of our algorithm makes the existing complexity bounds more precise for the case when the number of equations linking the positive definite variable matrices is not less than the sum of the ranks of the respective positive-semidefinite cones. 

\end{abstract}

\section{Introduction}

In the present paper, we propose a polynomial-time algorithm for convex feasibility problems with positive semidefinite constraints. Our algorithm is based on application of a family of projective transformations associated with certain valid inequalities, i.e., inequalities that are valid for solutions of the respective feasibility problem. 

The first application of projective transformations to linear programs refers to the well-known work \cite{karmarkar1984} of  Karmarkar. As noted later by Gill et al. \cite{gill1986} after the publication of Karmarkar's work, projected Newton barrier method and Karmarkar's algorithm are equivalent when the barrier parameter is suitably chosen. To the best of our knowledge, the previous use of projective transformations was restricted to linear optimization problems. Our goal in the present article is to show that an appropriate use of projective transformations leads to a polynomial algorithm in a more general context, namely, for convex feasibility problems with positive semidefinite constraints. Moreover, the running time of our algorithm makes the existing bounds more precise for the considered problem class.


According to the Koecher-Vinberg theorem, an element of a symmetric cone can be represented as a square of an element of the respective vector space in an associated Jordan algebra. The spectral theorem referring to Jordan algebras gives a theoretical basis for algorithmic solutions of convex feasibility problems with conic constraints involving symmetric cones. This makes it possible to develop a unified approach to semidefinite programming, second-order cone programming, and linear programming; see Lourenço et. al. \cite{lourenco2019}, Kanoh and Yoshise \cite{kanoh}, and Pe\~na and Soheili \cite{pena2016}.

The polynomial-time algorithms proposed in the mentioned works substantially depend on the Jordan-algebra constructions. In the present paper, we propose a simpler algebraic approach which can be applied to semidefinite programming problems stated in the form of linear systems with conic constraints. 

The main object of our further discussion is an abstract convex feasibility problem where the associated cone is assumed to be self-dual. We assume that its interior points induce linear operators having some additional properties. Then we define a class of projective transformations of a certain cross-section of the cone, which allows us to analyze the complexity of our algorithm. In the special case when the self-dual cone is the nonnegative orthant of a Euclidean space, the respective cross-section corresponds to the standard simplex and in the case when it is the cone of positive semidefinite (PSD) matrices, that cross-section corresponds to the spectraplex, which is formed by PSD matrices whose traces are equal to one. 

The algorithms proposed in \cite{kanoh}, \cite{lourenco2019}, and \cite{pena2016} for conic problems and the linear-programming algorithms proposed in \cite{chubanov2011} and \cite{chubanov2013} are based on two ingredients. One of them is what we will call a scaling procedure, which is a linear transformation of the problem in question, and a so-called basic procedure whose output is either a solution of the current problem or a valid inequality that holds for all feasible solutions in a ball with a given radius. This inequality is used in the framework of the scaling procedure to define the respective linear transformations which ensure the progress of the algorithm. The basic procedure is essentially the same in all of the mentioned works. It is closely related to Dantzig's algorithm for linear programming \cite{dantzig1992} and the relaxation method for linear inequalities by Agmon \cite{agmon1954} and Motzkin and Schoenberg \cite{motzkin1954} (which are exponential, see \cite{goffin1982}). 
In the present paper, we propose a scaling procedure which different from the previous works and is based on projective transformations. However, as an application we only consider semidefinite programming, whereas the mentioned results apply to abstract symmetric cones.

\section{Projective transformations of self-dual cones}\label{general-theory-section}

Let $K$ be a self-dual cone in a finite-dimensional real vector space $\cal H$ equipped with an inner product $\inprod{\cdot}{\cdot}.$ That is, \[x\in K\iff\inprod{x}{y}\ge 0, \forall y\in K.\]
Self-dual cones are necessarily full-dimensional, which means that the interior of $K,$ further denoted as $\interior{K},$ is not empty; see \cite{dattorro2005} for a reference.

Let $G$ be a class of linear operators from $\mathcal{H}$ to itself with the property that each operator in $G$ is uniquely identified by an element of $\interior{K}.$ If $a\in\interior{K}$ defines an operator of $G,$ we denote it by $L_a.$  
We assume that the following conditions are satisfied:
\begin{itemize}
\item[(i)] There exists an element $e\in\interior{K}$ (note that $\lambda e\in \interior{K}$ for any $\lambda > 0$) such that \[\forall a\in \interior{K}, \lambda > 0: L_a(e) = a,\;\;L_{\lambda e}(a) = \lambda a.\]
\item[(ii)] For each operator $g$ of $G:$
\begin{equation}\label{additional-property}
\forall a\in\mathcal{H}: \inprod{e}{g(a)} = \inprod{g(e)}{a}.
\end{equation}
\item[(iii)] For all $a\in \interior{K}:$ $e + a\in \interior{K},$ $L_{e + a}$ is invertible and $L^{-1}_{e + a}\in G.$ 
\item[(iv)] Each operator $g$ in $G$ maps $K$ into itself: $g(x)\in K$ for all $x\in K.$
\end{itemize}
\begin{lemma}\label{strict-ineq-observation}
 $\inprod{a}{b} > 0$ for all $a\in K,$ $a\neq {\bf 0},$ and $b\in \interior{K}$ 
\end{lemma}
{\bf Proof.} Since $b\in \interior{K},$ there is $\lambda > 0$ such that $b + \lambda(b - a) \in K.$ Then,
\[ 0\le\inprod{a}{b + \lambda(b - a)} = (1 + \lambda)\inprod{a}{b} -\lambda\|a\|^2,\] which implies that $\inprod{a}{b} > 0$ because $\lambda\|a\|^2 > 0.$
\BlackBox

\begin{lemma}\label{image-in-interior-lemma}
For all $g\in G:$ $g(e)\in\interior{K}.$
\end{lemma}
{\bf Proof.} If $g(e)\not\in\interior{K},$ then there exists $c\in K, c\ne{\bf 0},$ such that $\inprod{c}{g(e)} = 0,$ which implies, by (ii), $0 = \inprod{c}{g(e)} = \inprod{g(c)}{e}.$  Then Lemma \ref{strict-ineq-observation} implies that $g(c) = {\bf 0}$ because $g(c)\in K$ and $e\in\interior{K}.$ It follows that ${\bf 0} = g({\bf 0}) = g(c)$ while $c\ne{\bf 0},$ which means that $g$ is not invertible. We have a contradiction.  \BlackBox

For the purpose of analysis of the performance of our algorithm in the subsequent sections, below we define a class of projective transformations over the set
\[\Delta = \{x\in K: \inprod{e}{x} = 1\}.\]
which is a cross-section of $K.$
For instance, if $K = \mathbb{R}^n_+,$ then $\Delta$ is the standard simplex. If $K$ is the set $\psd{n}$ of $n\times n$ positive semidefinite (PSD) matrices, then $\Delta$ is a spectraplex. By this reason, we call $\Delta$ the spectraplex associated with $K,$ also in the general case.
For each element $y\in K$ we define the following projective transformation over $\Delta:$ 
\[
F_y(x) = \inprod{e}{L_{e + y}(x)}^{-1}L_{e + y}(x).
\]
We consider $F_y$ to be defined only on  $\Delta.$
As shown in the lemma below, $F_y$ is invertible and  the respective inverse transformation takes the form
\[
F^{-1}_y(x) = \inprod{e}{L^{-1}_{e + y}(x)}^{-1}L^{-1}_{e + y}(x).
\]

\begin{theorem}
$F_y(\Delta) = \Delta.$
\end{theorem}
{\bf Proof.} Let us show that \[\forall x\in \Delta: F_y(x)\in \Delta.\] Since $x, y\in K,$ it follows that $L_{e + y}(x)\in K.$ Then $F_y(x) \in K$ because $F_y(x) = \lambda L_{e + y}(x)$ where $\lambda = (1 + \inprod{y}{x})^{-1} > 0.$ ($\inprod{y}{x} \ge 0$ because $K$ is self-dual). At the same time, $\inprod{e}{F_y(x)} = 1.$

Now we consider $x^\prime\in\Delta$ and show that there exists $x\in\Delta$ with $F_y(x) = x^\prime.$ Indeed, let $x = F_y^{-1}(x^\prime).$ 
Lemma \ref{strict-ineq-observation} implies that \[\inprod{e}{L^{-1}_{e + y}(x^\prime)} = \inprod{L^{-1}_{e + y}(e)}{x^\prime} > 0\] because, by Lemma \ref{image-in-interior-lemma}, $L^{-1}_{e + y}(e)\in \interior{K}$ and  $x^\prime\in K\setminus\{\bf 0\}.$ Then $x\in K$ because $L^{-1}_{e + y}(x^\prime) \in K.$ At the same time, $\inprod{e}{x} = 1.$ It follows that $x\in\Delta.$ \BlackBox

Let $A$ be a linear operator and ${\bf 0}$ denote a zero in the codomain of $A.$ (Further, whenever we use notation ${\bf 0},$ it should be clear from the context what space it refers to.)
A problem of the form
\begin{equation}\label{problem-in-question}
Ax = {\bf 0}, x\in \interior{K},
\end{equation} has a feasible solution if and only if the problem
\[Ax = {\bf 0}, \inprod{e}{x} = 1, x\in \interior{K}\]
has a feasible solution. To obtain a solution of this problem from a solution $x$ of problem (\ref{problem-in-question}), we can simply divide it by the value $\inprod{e}{x},$ which must be positive by Lemma \ref{strict-ineq-observation}. That is, if (\ref{problem-in-question}) is feasible, then it has a solution in $\Delta.$ 


Observe that 
\[
\inprod{e}{L_{e + y}(x)} = \inprod{L_{e + y}(e)}{x} = \inprod{e + y}{x} = \inprod{e}{x} + \inprod{y}{x} = 1 + \inprod{y}{x}.
\]
It follows that
\[
F_y(x) = \frac{L_{e + y}(x)}{1 + \inprod{y}{x}}.
\]
That is, the smaller the value of $\inprod{y}{x},$ the closer $F_y(x)$ is to $L_{e + y}(x).$
A key ingredient of our algorithm is a procedure which either finds a feasible solution to problem (\ref{problem-in-question}) or a sufficiently small $\varepsilon$ such that
\begin{equation}\label{valid-inequality}
\forall x\in X^*\cap\Delta: \inprod{y}{x}\le\varepsilon,
\end{equation}
where $X^*$ is the solution set of (\ref{problem-in-question}).

If an inequality of the above form is found, our algorithm transforms the problem by means of 
replacing $A$ by $A^\prime = AL_{e + y}^{-1}.$ The new problem has the form \[A^\prime x = {\bf 0}, x\in\interior{K}.\]  
The solution set of the new problem is $L_{e + y}(X^*),$ which follows from (iv). The described transformation of the problem is called a {\em scaling iteration}. To go back from a solution $x^\prime$ of the new problem to a solution of the original problem, it suffices to apply $L^{-1}_{e + y}$ to $x^\prime.$

To measure the progress of our algorithm, we observe what happens with solutions in $\Delta.$ If $x\in\Delta $ is a solution of the original problem (\ref{problem-in-question}), then $x^\prime = F_y(x)$ is a solution of the new problem, in the same set $\Delta.$
In the next section we introduce the so-called potential function $\varphi$ bounded on $\Delta.$ (We call $\varphi$ the potential function by analogy with the potential function used by Karmarkar to measure the progress of his projective method in each iteration); when $L_{e + y}$ is applied to $x\in X^*\cap\Delta,$ the value of $\varphi$ grows by a constant factor depending only on the dimension of the space. Then only a small correction is needed (division by $1 + \inprod{y}{x}$) to bring the resulting vector back to $\Delta,$ which in turn requires only a small correction of that constant factor to estimate the increment of $\varphi$ when applying the projective transformation $F_y$ to $x.$ Since $\varphi$ is bounded on $\Delta,$ this provides a tool for estimating the number of scaling iterations of our algorithm. 

%
\section{Potential function}

Let $\varphi: {K}\longrightarrow\mathbb{R}_{+}$ be a continuous function 
such that \[\forall a\in\interior{K}: \varphi(a) > 0,\]
\begin{equation}\label{phi-inequality}
\forall a, b\in \interior{K}: \varphi(L_a(b)) \ge \varphi(a)\varphi(b),
\end{equation}
and 
\begin{equation}\label{range-inequality}
\varphi({\bf 0}) = 0, \varphi(e) = 1.
\end{equation}
We also assume the existence of $\eta_\varphi > 0$ such that
\[
\inf_{y\in\Delta} \varphi(e +y) \ge 1 + \eta_\varphi.
\]
In the next sections we provide an explicit formula for $\eta_\varphi$ for the case of semidefinite programming.


\begin{lemma}\label{phi-progress-observation} For all $x,y\in\Delta:$
\begin{equation}\label{progress-ineq}
\varphi(F_y(x)) \ge \varphi(e + y)\cdot\varphi((1 + \inprod{x}{y})^{-1}e)\cdot\varphi(x).
\end{equation}
\end{lemma}
{\bf Proof.} 
Observe that $\lambda a = \lambda L_a(e) = L_a(\lambda e).$ Note that 
\[
\forall x\in\Delta: \inprod{e}{L_{e + y}(x)} = \inprod{L_{e + y}(e)}{x} =  \inprod{e + y}{x} = 1 + \inprod{y}{x}.
\]
Then, using (i), we obtain 
\[
F_y(x) = L_{e + y}((1 + \inprod{x}{y})^{-1}x) = L_{e + y} L_{(1 + \inprod{x}{y})^{-1}e}(x).
\]
Now, applying (\ref{phi-inequality}) two times, we get (\ref{progress-ineq}).
\BlackBox

Lemma \ref{phi-progress-observation} implies that 
\[
\varphi(F_y(x)) \ge \varphi(x)\cdot(1 + \eta_\varphi)\cdot\varphi((1 + \inprod{x}{y})^{-1}e).
\]
Since $\varphi$ is continuous, it follows from (\ref{range-inequality}) that for each $\kappa\in ]0,1[$ there exists $\varepsilon$ such that 
\begin{equation}\label{epsilon-kappa}
\forall x,y\in \Delta \mbox{ with } \inprod{y}{x} \le \varepsilon: \varphi((1 + \inprod{y}{x})^{-1}e) \ge 1 - \kappa.
\end{equation}
We set
\begin{equation}\label{kappa}
\kappa := \frac{\eta_\varphi}{2\cdot(1 + \eta_\varphi)}.
\end{equation}
Then,
\[
(1 + \eta_\varphi)\cdot (1 - \kappa) = 1 + \eta_\varphi - (1 + \eta_\varphi)\cdot\kappa = 1 + \frac{\eta_\varphi}{2}.
\]
It follows that 
\begin{equation}\label{progress-inequality}
\forall x,y\in \Delta \mbox{ with } \inprod{y}{x} \le \varepsilon: \varphi(F_y(x)) \ge \varphi(x)\left(1 + \frac{\eta_\varphi}{2}\right).
\end{equation}
Further, we assume that $\kappa$ is chosen exactly as in (\ref{kappa}) and the $\varepsilon$ is chosen so as to ensure (\ref{epsilon-kappa}). Our algorithm is based on a procedure which either finds a feasible solution or a suitable $y\in\Delta$ such that (\ref{epsilon-kappa}) holds for all feasible $x.$
\section{Abstract algorithm}

Since $K$ is self-dual, problem (\ref{problem-in-question}) can equivalently be written as
\begin{equation}\label{problem}
Ax = \zero, \inprod{u}{x} > 0, \forall u\in \Delta.
\end{equation}
We assume that $\Delta$ is bounded.
Let $\theta(\Delta) < \infty$ be an upper bound for $\|x\|, x\in \Delta.$ In the case when $\Delta$ is a standard simplex or the set of all positive semidefinite matrices with trace equal to $1,$ one can simply set $\theta(\Delta) = 1.$ 

Let $P_A$ denote the projection operator onto $\ker A$ (the nullspace of $A$).  Assume that we are given an upper bound $U^+$ for $\varphi$ over $\Delta,$ and a lower bound $U^-$ that must hold for $\varphi(x^*)$ at some feasible solution $x^*\in\Delta$ if the problem is feasible. 
The above observations and our previous discussion on $\varphi$ lead to the following abstract algorithm:

\begin{algorithm}\label{abstract-algorithm}
\end{algorithm}
\begin{algorithmic}
\STATE{{\bf Input:} Linear operator $A.$}
\STATE{{\bf Output:} A solution to (\ref{problem}) or a decision that (\ref{problem}) is infeasible.}
\STATE{$A^\prime := A,$ $k := 0.$}
\STATE{Initialize $y$ with an arbitrary element of $\Delta$ and $L$ with the identity map.}
\WHILE{$(1 + \eta_\varphi/2)^k < U^+/U^-$}
\IF{$\inprod{u}{P_{A^\prime}y}  > 0$ for all $u\in\Delta$}
\STATE{Return $L^{-1}P_{A^\prime}y.$ (It is feasible for (\ref{problem}) and (\ref{problem-in-question}).)}
\ELSE
\STATE{Pick $u\in\Delta$ with $\inprod{u}{P_{A^\prime} y} \le 0.$}
\ENDIF
\STATE{$\alpha := \inprod{P_{A^\prime}u}{P_{A^\prime}(u - y)}/\|P_{A^\prime}(u - y)\|^2$}
\STATE{$y := \alpha y + (1 - \alpha)u$}
\IF{$\|P_{A^\prime} y\| \le \varepsilon/\theta(\Delta)^2$}
\STATE{$L := L_{e + y}L,$ $A^\prime := A^\prime L^{-1},$ $k := k + 1.$}
\ENDIF
\ENDWHILE
\STATE{Return that (\ref{problem}) has no solutions.}
\end{algorithmic}
In the algorithm, the $\alpha$ is chosen so that $\alpha P_{A^\prime}y + (1 - \alpha)P_{A^\prime}u$ is the orthogonal projection of $\zero$ on the segment $[P_{A^\prime}y, P_{A^\prime}u].$ Note that $\inprod{P_{A^\prime}y}{P_{A^\prime}u} = \inprod{P_{A^\prime}y}{u}\le 0,$ which follows from the fact that $P_{A^\prime}$ is self-adjoint and idempotent. Therefore, the respective step guarantees that $\|P_{A^\prime}y\|^{-2}$ increases by at least a specified minimum value, as follows from the reciprocal Pythagorean theorem: 
\begin{lemma}\label{reciprocal-corollary}[Corollary of the reciprocal Pythagorean theorem]
Let $a$ and $b$ be vectors of an inner-product space such that $\inprod{a}{b} \le 0.$ Let $h$ be the  projection of $\zero$ on the segment $[a,b].$ Then,
\[\norm{h}^{-2} \ge \norm{a}^{-2} + \norm{b}^{-2}.\]
\end{lemma}
{\bf Proof.} Using $\inprod{a}{b} \le 0,$ by a direct calculation we can make sure that $\inprod{a - b}{h} = 0,$ which means that $\|h\|$ can be viewed as the height of the triangle formed by $a$ and $b.$ Let $c_1 = a - h$ and $c_2 = b - h.$ Then, by the Pythagorean theorem,

\[\norm{h}^2 = \norm{a}^2 - \norm{c_1}^2,\;\;\;\norm{h}^2 = \norm{b}^2 - \norm{c_2}^2,\;\;\;(\norm{c_1} + \norm{c_2})^2 = \norm{a}^2 + \norm{b}^2 - 2\inprod{a}{b}.\]
Combining these equations, we see that
\[2\norm{h}^2 = \norm{a}^2 + \norm{b}^2 - \norm{c_1}^2 - \norm{c_2}^2 = 2\inprod{a}{b} + 2\norm{c_1}\norm{c_2} \le 2\|c_1\|\|c_2\|.\]
It follows that
$\norm{h}^4 \le \norm{a}^2\norm{b}^2 - \norm{h}^2(\norm{a}^2 + \norm{b}^2) + \norm{h}^4.$

Then,
\[\norm{h}^2 \le \frac{\norm{a}^2\norm{b}^2}{\norm{a}^2+\norm{b}^2} \Longrightarrow \norm{h}^{-2} \ge \norm{a}^{-2} + \norm{b}^{-2}.\]
\BlackBox

By an iteration of Algorithm \ref{abstract-algorithm}, we mean an iteration of its while-loop. An iteration in which the condition $\|P_{A^\prime}y\| \le \varepsilon/\theta(\Delta)^2$ holds is called a {\bf scaling iteration}.
\begin{lemma}\label{basic-procedure-complexity-lemma} 
Algorithm \ref{abstract-algorithm} performs no more than $O(\theta(\Delta)^4/\varepsilon^2)$ iterations between two consecutive scaling iterations.
\end{lemma}
{\bf Proof.} 
Consider an iteration of the algorithm. Let $y^\prime$ be the new vector $y$ and $y$ itself be the old $y.$ The $\alpha$ is chosen so as to make $P_{A^\prime} y^\prime$ the projection on $[P_{A^\prime} y, P_{A^\prime} u].$ Then, Lemma \ref{reciprocal-corollary} implies that
\[\|P_{A^\prime} y^\prime\|^{-2}\ge  \|P_{A^\prime} y\|^{-2} + \|P_{A^\prime} u\|^{-2} \ge \|P_{A^\prime} y\|^{-2} + \|u\|^{-2} \ge \|P_{A^\prime}y\|^{-2} + \theta(\Delta)^{-2}.\] 

After $l$ iterations from the previous to the next scaling, we get
$
\|P_{A^\prime} y^\prime\|^{-2} \ge l\theta(\Delta)^{-2}.
$
The next scaling takes place no later than when  $l \ge \theta(\Delta)^4/\varepsilon^{2}.$
\BlackBox

\begin{theorem} Algorithm \ref{abstract-algorithm} either finds a feasible solution or decides that the original problem is infeasible in \[O(\varepsilon^{-2}\theta(\Delta)^4\log(1 + \eta_\varphi/2)^{-1}\log (U^+/U^-))\] iterations.
\end{theorem}
{\bf Proof.} Consider an iteration of Algorithm \ref{abstract-algorithm}. Let $x\in\Delta$ be a feasible solution of the current problem, i.e, $x\in \ker{A^\prime}\cap\Delta.$ If $\|P_{A^\prime}y\| \le \varepsilon/\theta(\Delta)^2,$ then \[\inprod{y}{x} = \inprod{y}{P_{A^\prime}x} = \inprod{P_{A^\prime}y}{x} \le \|P_{A^\prime}y\|\|x\| \le \varepsilon.\]
Then $x$ and $y$ satisfy (\ref{progress-inequality}). 
Now let $x^*$ be a feasible solution of the original problem with $\varphi(x^*) \ge U^-.$ Let $x^*_k\in\Delta$ be a feasible solution obtained from $x^*$ by applying the respective sequence of projective transformations defined by $y_i\in \Delta,$ $i\in[k],$ corresponding to the $y$'s used to define linear operators $L_{e + y}$ at the respective scaling iterations. Thus, after $k$ scaling iterations we have 
\[
U^+\ge\varphi(x^*_k) \ge \varphi(x^*)\left(1 + \eta_\varphi/2\right)^k \ge U^-\left(1 + \eta_\varphi/2\right)^k.
\]
It follows that
\[
k \le \left\lceil\log(1 + \eta_\varphi/2)^{-1}\log (U^+/U^-)\right\rceil,
\]
which is an upper bound on the number of scaling iterations. The number of iterations between two scaling iterations is bounded by $O(\theta(\Delta)^4/\varepsilon^2),$ which follows from Lemma \ref{basic-procedure-complexity-lemma}. 
 \BlackBox


\section{Application to semidefinite programming}

Now we consider $\mathcal{H} = V_1\times\ldots\times V_l,$ where each $V_s$ is the space $\mathbb{S}^{n_s}$ of symmetric matrices. In this section, the abstract problem (\ref{problem-in-question}) takes the more specific form of a feasibility problem with $m$ linear equations linking positive definite variable matrices:
\[
\sum_{s\in[r]} A_sx_s = {\bf 0}, x_s\in \interior{K_s}, s\in[l],
\]
where $A_s:V_s\map\mathbb{R}^m$ is a linear operator and $K_s,$ $s\in[l],$ is the respective cone $\psd{n_s}$ of positive semidefinite matrices. 

Let $K = K_1\times\ldots \times K_l.$ 
Cone $K$ is self-dual because $K_s$ are self-dual. Note that $n_s$ can be equal to $1,$ in which case $K_s = \mathbb{R}_+$ and $\interior{K_s} = \mathbb{R}_{++}.$

Thus, vectors $x,y\in \mathcal{H},$ are tuples $x = (x_1,\ldots,x_l)$ and $y = (y_1,\ldots,y_l)$ whose components are square matrices in the respective spaces. Let the product $x\circ y$ be defined as $x\circ y = (x_1\cdot y_1,\ldots, x_l\cdot y_l),$ where $\cdot$ is the usual matrix product.

An eigenvalue of $x\in \mathcal{H}$ is simply an eigenvalue of one of its components $x_s.$ The vector of the eigenvalues of $x$ is composed of the eigenvalues of its components.


We define the trace of $x$ as the sum of traces of its components and the determinant as the product of determinants:
\[\tr{x} := \sum_{s\in[l]}\tr{x_s},\;\;\; \det{x} := \prod_{s\in[l]}\det{x_s},\]
where $\tr{x_s}$ is the trace of matrix $x_s$ and $\det{x_s}$ is its determinant. 

The well-known cyclic permutation property of traces of matrices implies that 
\[\tr{(x\circ y\circ z)} = \tr{(y\circ z\circ x)}.\]
The inner product $\inprod{x}{y}$ of $x$ and $y$ is naturally defined as 
\[\inprod{x}{y}:=\tr{(x\circ y)}.\] 
The inner product induces the norm $\|x\| = \inprod{x}{x}^\onesecond.$


Since PSD matrices can be represented as squares of PSD matrices, an $
y\in K$ can be factored as $y = y^\onesecond\circ y^\onesecond,$ where $y^\onesecond\in K.$

Note that $y\in\interior{K}$ is invertible w.r.t. operation $\circ.$ Indeed, $y^{-1}$ in that case is simply $y^{-1} = (y^{-1}_1,\ldots, y^{-1}_l)$ where each $y^{-1}_s$ is the respective inverse matrix. 

Let the following linear operator $L_a$ be associated with $a\in\interior{K}:$
\begin{equation}\label{semidefinite-operation}
L_a(x) = a^\onesecond\circ x\circ a^\onesecond.
\end{equation}
Note that $L^{-1}_a = L_{a^{-1}}.$
\begin{lemma}
Linear operators defined by (\ref{semidefinite-operation}) satisfy (i)-(iv).
\end{lemma}
{\bf Proof.} 
Conditions (i)-(iv) are implied by the following observations:

(i): The element $e$ can be defined as $e = (I_1,\ldots,I_r),$ where $I_s$ is the identity matrix in $V_s.$ Then, $L_a(e) = a$ for all $a\in\interior{K}.$ 


(ii): Observe that \[\inprod{e}{L_y(x)} = \tr{(y^\onesecond\circ x\circ y^\onesecond)} = \tr{(x\circ y^\onesecond\circ y^\onesecond)} = \tr{(x\circ y)} = \inprod{y}{x} = \inprod{L_y(e)}{x},\] where we have applied the mentioned cyclic permutation property of the trace. 

(iii): This follows from the invertibility of $e + a$ for any $a\in K$ and the fact that $L^{-1}_{e + a} = L_{(e + a)^{-1}}.$

(iv): This property follows directly from the fact that the matrix components of $a^\onesecond \circ x \circ a^\onesecond$ are PSD matrices.
%
 \BlackBox

%
%
%
%

The set $\Delta$ takes the following form:
\[
\Delta = \{x\in K: \inprod{e}{x} = 1\} = \{x\in K: \tr{x} = 1\}.
\]
In the case $l = 1,$ this set is a spectraplex whose center is the identity matrix divided by $n.$ It is natural to call $\Delta$ the spectraplex also in the case $l > 1,$ with the center $e/n.$

For any $x\in\Delta,$ we have $\|x\| \le \tr{x} = 1.$ Therefore, we set $\theta(\Delta) = 1.$ 

The respective projective transformation induced by $y\in K$ is
\[
F_y(x) = \inprod{e}{L_{e + y}(x)}^{-1} L_{e + y}(x) = (1 + \inprod{y}{x})^{-1} (e + y)^\onesecond\circ x\circ (e + y)^\onesecond.
\]

We define the potential function $\varphi$ as the determinant:
\[
\varphi(x) := \det{x}.
\]
Note that the $e/n$ delivers the unique maximum of $\varphi$ over $\Delta.$ So an element $x$ of $\Delta$ with a higher value of $\varphi(x)$ lies closer to $e/n$ in the sense of a similarity measured in terms of the difference between $\varphi(e/n) - \varphi(x).$
All the required properties follow from the properties of determinants of matrices, in particular, 
\[\det{L_{e + y}(x)} = \det{((e + y)^\onesecond\circ x \circ (e + y)^\onesecond)} = \det{(e + y)}\cdot\det{x}.\]
For $y\in\Delta,$ the eigenvalues $\lambda_1,\ldots,\lambda_n$ of $y$ are nonnegative and \[1 = \inprod{e}{y} = \sum^n_{j = 1}\lambda_j.\]
Therefore, one can choose $\eta_\varphi = 1$ because \[\det{(e + y)} = \prod_{j = 1}^n (1 + \lambda_j) \ge 1 + \sum^n_{j = 1}\lambda_j = 2.\]
Then, following the formula (\ref{kappa}) for $\kappa,$ we calculate \[\kappa = 1/4.\] 
Let 
\[ n = \sum_{s\in[l]} n_s.\]
Now we choose $\varepsilon$ as \[\varepsilon := \frac{\ln(4/3)}{n}\]
and make sure that this choice of $\varepsilon$ is suitable for (\ref{epsilon-kappa}):
\[
\det{((1 + \inprod{y}{x})^{-1} e)} = (1 + \inprod{y}{x})^{-n} \ge (1 + \varepsilon)^{-n} \ge 3/4 = 1 - \kappa.
\]
Then (\ref{progress-inequality}) takes the form
\[
\varphi(F_y(x)) \ge \varphi(x)\left(1 + \frac{\eta_\varphi}{2}\right) = \frac{3}{2}\varphi(x).
\]

Assume that there is a feasible solution $x^*\in X^*\cap\Delta$ and $\lambda_{\min}(x^*)$ is its minimum eigenvalue. Then we can choose $U^- = \lambda_{\min}(x^*)^n$ as a positive lower bound for $\det{x^*}.$ Since $\det x \le \varphi(e/n)$ for all $x\in\Delta,$ we set $U^+ = 1/n^n.$ 

Now we apply Algorithm \ref{abstract-algorithm}. The number of scaling iterations is bounded by
$
O(n\log_{3/2} 1/(n\lambda_{\min}(x^*))).
$
The number of consecutive iterations between two scaling iterations is bounded by $O(n^2).$
Besides taking the square root of $e + y,$ the most computationally expensive operations are: The computation of the projection on $\ker{A^\prime}$ and finding a suitable $u\in\Delta$ in the course of the algorithm. 
The $u$ can be computed via an eigenvector $v$ of $P_{A^\prime}y$ corresponding to its minimum eigenvalue ($u := vv^T/\|v\|^2$), in the same way as proposed in \cite{lourenco2019}. 

With respect to the running time, the worst-case scenario is when $l = 1,$ in which case $K = \mathbb{S}^n_+.$
As well as in \cite{lourenco2019}, we consider a model of computation where the minimum eigenvalue of a symmetric matrix can be found in time $c_{\min}.$ Then a suitable $u$ in the course of the algorithm can be found in time $O(n^{2.5} + c_{\min})$ using fast matrix inverse. For a given symmetric matrix its eigenvalue decomposition with pairwise orthogonal eigenvectors can be found in $O(nc_{\min} + n^3)$ time.

In the above estimate we have taken into account the following computations:
\begin{itemize}
\item The computation of the projection operator $P_{A^\prime}:$ This is needed only at scaling iterations and each time requires $O(n^2m^2)$ time.
\item As mentioned above, the computation of $u$ takes $O(n^{2.5} + c_{\min})$ time, i.e., an overall time of $O(n^{4.5} + n^2c_{\min})$ between two consecutive scaling iterations.
\item Taking the square root of $e + y$ at each of the scaling iterations requires $O(n^3 + nc_{\min})$ time.
\end{itemize}

{\em Remark.} Note that $m$ can be comparable with or greater than $n.$ This is, for instance, a typical situation the case when a semidefinite program serves as a relaxation of a combinatorial problem.

More precisely, in the worst-case scenario with $l = 1,$ our algorithm runs in time 
 \[O\left(\inf_{x^*\in X^*\cap\Delta}n\log\frac{1}{n\lambda_{\min}(x^*)}\cdot(m^2n^2 + n^{4.5} + n^2c_{\min})\right)\]
 in the case when the solution set $X^*$ is nonempty.

A typical task associated with the feasibility problem we consider in this section is to either find a solution or prove that no solution exists with the minimum eigenvalue less than a given positive value $\delta.$ In fact, we can approximate our problem with a given accuracy so as to ensure existence of such a solution in case of feasibility. For this, we can replace the operator equation $Ax = {\bf 0}$ by $A\bar{x} - \delta Ae = {\bf 0},$ to guarantee that the resulting approximate problem has a solution with the minimum eigenvalue not less than $\delta/(n + 1)$ in $\Delta.$ The homogenization leads to the problem $A\bar{x} - \delta\bar{x}_{r + 1}Ae = {\bf 0}, \bar{x}\in\interior{K}, \bar{x}_{r+1} > 0,$ to which we can apply our algorithm; the cone $K$ is replaced by $K\times\mathbb{R}_{++}$ in this case.

So, for the above task, our algorithm runs in time
 \[O\left(n\log(1/\delta)\cdot(m^2n^2 + n^{4.5} + n^2c_{\min})\right)\] when no solutions in $\interior{K}$ exist. 
For the same task, the running time of the algorithm proposed by  Lourenço at al. in \cite{lourenco2019} is bounded by \[O(n\log(1/\delta)(m^2n^2 + \max\{mn^4, n^2c_{\min}\}).\]
For the case when $m \ge n,$ the running time of our algorithm is slightly better.

In the general case, when $K$ is a Cartesian product of $l$ PSD cones, our algorithm runs in time
 \[O\left(n\log(\delta^{-1})\cdot\left (m^2\sum_{s\in[l]}n^2_s + \sum_{s\in [l]}n^{4.5}_s + n^2c_{\min}\right)\right)\]
In this case, the algorithm proposed in 
\cite{lourenco2019} has a running time of
\[
O\left(\frac{r}{\gamma}\log\left(\delta^{-1}\right)(m^3 + m^2d + l^3r^2_{\max}\cdot\max\{md, c_{\min}\})\right),
\]
where $r$ is the rank of the symmetric cone, $r_{\max}$ is the maximum rank of a simple symmetric cone participating in the Cartesian product of $l$ simple symmetric cones representing the original one, $\gamma$ is a value not greater than $2,$ and $c_{\min}$ is the computation cost of finding the minimum eigenvalue of an element of a symmetric cone.  For our problem, one should set $r = n$ and $r_{\max} = \max_{s\in[l]}n_s,$ and $d = n(n-1)/2 + n.$ 


\section{Summary}

We have studied a class of projective transformations and used it for the complexity analysis of our algorithm for an abstract convex feasibility problem over a self-dual cone. On this basis, we have proposed a polynomial algorithm for feasibility problems with positive definite constraints. In the case when the number of equations linking positive definite variable matrices is greater than the sum of their dimensions, we have slightly improved the existing running time. 

\end{document}